\documentclass[11pt]{article}

\usepackage{amssymb,amsmath,latexsym,epsfig,euscript}

\def \N{\mathbb{N}}

\def \R{\mathbb{R}}
\def \Z{\mathbb{Z}}

\newtheorem{defn}{Definition}[section]

\newtheorem{prop}[defn]{Proposition}

\newenvironment{undef}[1]%
           {\vspace{3.3mm}
           \noindent{\bf #1}\it}%
           {\vspace{3.3mm}}

\newenvironment{proof}[1]{
  \trivlist \item[\hskip \labelsep{\it #1}]}{\hfill\mbox{$\square$}
  \endtrivlist}

\begin{document}

\hfill May 2nd., 2001

\vspace{4mm}

\noindent{\LARGE{\bf Intrinsic palindromic numbers}}

\vspace{7mm}

\noindent {\large Antonio J. Di Scala and Mart{\'\i}n Sombra }


\vspace{5mm}


\noindent {\small {\bf Abstract.}
We introduce a notion of palindromicity of a
natural number which is independent of the base.
We study the existence and density of palindromic and multiple palindromic
numbers, and we raise several related questions.
}

\vspace{1mm}


\vspace{1mm}


\vspace{-8mm}


\setcounter{section}{1}

\section*{}

Natural numbers occur everywhere in our daily life: a bus
ticket, a car plate, an id number, a timetable, etc..
These numbers are mostly expressed in the decimal system.
That is,  for a natural number $n \in \N$
we write
$n = (a_{k-1} \, \cdots
a_0)_{10}$
for some $0 \le a_i \le 9$ and $a_{k-1} \ne 0$, which means that
$$
n = a_{k-1} \, 10^{k-1}+ a_{k-2} \, 10^{k-2}+ \cdots + a_0 .
$$

\smallskip

Of particular attraction are the so-called palindromic numbers.
These are the numbers whose decimal expansion is the same when
read from left to right, and from right to left, that is
$(a_{k-1} \, \cdots
a_0)_{10}  = (a_0 \, \cdots
a_{k-1})_{10}$.

This kind of numbers appears already in the
Ganitas{\^a}rasamgraha, a sanskrit manuscript
dated around 850 AD.
Therein the Indian mathematician Mah{\^a}v{\^\i} r{\^a}ch{\^a}rya
described the 12345654321 as the quantity which ``beginning with
one until it reaches six, then decreases in reverse order''
\cite[p.399]{Ifrah00}.
This is a curious palindromic number, and in particular it is the square of another
palindrome:
$12345654321 = 111111^2$.
There are many ways of generating palindromic numbers
(see for instance \cite{Reichmann57}) including an interesting conjecture
\cite{Trigg67}.

\bigskip

Crossing by chance with one such a number
is a rare occurrence: for instance, the
probability of
picking at random a number  $ 10^4 \le n < 10^5$
and it resulting a palindrome is $1/10^2$, a fact
well-known by the collectors of palindromic bus
tickets.
Hence we tend to feel pretty lucky when one of this rare numbers
crosses our path.
And the more figures it has, the luckier we
feel, and the luckier the number itself seems to be.

\smallskip

But, is this feeling really justified?
The truth has to be said:
this property is {\em not } intrinsic  to the
number, but also depends on the base used to express it.
So the number $n:= (894111498)_{10}$ is lucky (or palindromic) in base
10, but is not in base 13 as
$
n=(113314377)_{13} $.
Hence each time we encounter the 894111498  we have to --- at least in principle ---
thank heaven for this
occurrence {\em together } with the fact that the human race has five
fingers in each hand.

\smallskip

However we can easily get independent of the base,
defining a number to be intrinsically
palindromic if is palindromic
in {\em some } base.
A minute reflection shows that this definition
is meaningless as it stands:
every number is palindromic in any base $m >n$, as  $n=(n)_m$.

Indeed it is much more natural to take into account the number of
figures.
So we define a number $n \in \N$ to be
{\em $k$-palindromic } if there is a base $b$ such that the $b$-expansion of
$n$ is palindromic of length $k$.

The previous observation shows that every number is
1-palindromic.
Also note that for all $n\ge 2$ it holds $n=(1 1)_{n-1}$,  that is
every $n\ge 2$ is 2-palindromic.

\medskip

What about $k$-palindromic numbers for $k \ge 3$?
Our thesis is that very few numbers are
$k$-palindromic, at least for  $k \ge 4$:
the probability of a number in the appropriate range
being, say, intrinsically
$9$-palindromic is small, and indeed quite close to the
probability of being
$9$-palindromic in base 10.
This justifies our first impression that the 894111498 is lucky,
regardless the base chosen to represent it.

\smallskip

To write down our results we first have to introduce the following
counting functions.
Take $k, N , b \in \N$ and set
$$
\Phi_k(N, b):= \# \{ \ n \le N \, ; \,  \ n \mbox{ is } k\mbox{-palindromic in
base } b \  \}.
$$
Also set
$$
\Phi_k(N):=
\# \{ \ n \le N \, ; \,  \ n \mbox{ is } k\mbox{-palindromic} \  \}.
$$
Then $\Phi_k(N,b) / N$ and $\Phi_k(N) / N$ stand for
the {\em density } (or probability) of numbers below
$N$
which are $k$-palindromic
in base $b$ and
intrinsically $k$-palindromic,
respectively.

\begin{undef}{Theorem 1} \label{Phi_k}
Let $k \ge 4$,
and write $k = 2 \, i + r$ with $ i \in \N$ and $r =0, 1$.
Then
$$
\Phi_k(N) \le 4 \, (N+1)^{\frac{i+r  +1}{k}}.
$$
\end{undef}

\vspace{-11mm}

\begin{proof}{Proof.--}
A base $b$ contributes to $\Phi_k (N)$ if and only if $\Phi_k(N,b)
>0$, namely if and only if there exists a number
$$
n=
a_{k-1} \, b^{k-1} + \cdots +
a_0 \ \le N $$
palindromic in base $b$ of length
$k$, that is such that $ 0 \le a_j \le b-1$ and
$a_{k-j } = a_j$ for $j = 0, \dots, i+r $,  and
$a_{k-1} \ne 0$.
Then  $b^{k-1} +1 \le n \le N$, and hence  $b$ contributes to
$\Phi_k(N)$ if and only if
$ b \le (N-1)^{\frac{1}{k-1}}$.

\smallskip

We consider separately the cases $b \le (N+1)^{\frac{1}{k}}$ and $ (N+1)^{\frac{1}{k}}
< b \le (N-1)^{\frac{1}{k-1}}$.
In the first case, the largest
$k$-palindrome in base $b$ is
$$
(b-1) \, (b^{k-1} + \cdots + 1)= b^k -1 \le N,
$$
and so
$ \Phi_k(N,b)  = \Phi_k( \infty ,b) = (b-1)
\, b^{  i +r -1}  $.

\smallskip

Now for the second case we let $\theta(b) \in \N$ be the largest integer such that
$
\theta(b) \, (b^{k-1}+1) \le N$,
that is
$\theta(b) = [ N / (b^{k-1}+1) ]  \le  N / b^{k-1}$.
Then $(\theta(b)+1) \, (b^{k-1}+1) > N $ and so
every
$k$-palindromic number in base $b$ begins with $a_{k-1} \le
\theta(b)$.
Hence
$$
\Phi_k(N, b)  \le \theta(b) \, b^{i+r -1}  \, \le N / b^{i } .
$$
\vspace{-3mm}
Set
$$
L:= [(N+1)^{\frac{1}{k}}]  \quad \quad , \quad \quad M:= [(N-1)^{\frac{1}{k-1}}].
$$
\smallskip
We split the sum
$
\Phi_k(N) = \varphi + \chi + \psi $ with $\varphi :=
\sum_{b=2}^{L-1}
\Phi_k(N,b)$,
$\ \psi := \sum_{b=L+2 }^{M}  \Phi_k(N,b)$, and
$\ \chi := \Phi_k(N,L) + \Phi_k(N,L+1)$.
From the previous considerations we deduce that
$$
\varphi \, \le \, \sum_{b=2}^{L-1 }  b^{i+r}
\, \le \,
\int_2^L t^{i+r  } \, dt
\, \le \,
\frac{L^{i+r+1}}{i+r +1}
\, \le \,
\frac{(N+1)^{\frac{i+r+1}{k}}}{i+r +1}.
$$

On the other hand
$$
\psi
\, \le \,
\sum_{L+2}^{M} \frac{N}{b^i }
\, \le \,
N \, \int_{L+1 }^{M  } \frac{dt}{t^i }
\, \le \,
\frac{N}{(i-1) \, (L+1)^{i-1}}
\, \le \,
\frac{(N+1)^{\frac{i+r+1}{k}} }{i-1} .
$$

Finally $\, \chi \le L^{i+r} + N / (L+1)^i \le 2 \, (N+1)^\frac{i+r}{k}$
and thus
$$
\Phi_k(N)
\, \le \,
\left(
\frac{2}{(N+1)^\frac{1}{k}}+ \frac{1}{i+r+1} + \frac{1}{i-1}
\right) \,
(N+1)^{\frac{i+r+1}{k}}
\, \le \,
4 \,
(N+1)^{\frac{i+r+1}{k}} .
$$

\end{proof}

Let $k= 2\, i + r \ge 4$, and let $b\ge 2$ be a base. Set $N:=
b^k-1$, so that
$N$ is larger than every number whose representation in base $b$
has length $k$.
Then
$$
\Phi_k(N,b) =
\Phi_k(\infty,b) = (b-1) \, b^{i+r-1},
$$ and so the density of numbers below $N$
which are $k$-palindromic in base $b$ is $(b-1) \, b^{i+r-1} / N
\sim 1/b^i$.
On the other hand, the previous result shows that the density of
intrinsic $k$-palindromes below $N$
is bounded by $ 4 \, (b^k)^{\frac{i+r  +1}{k}}/ (b^k-1)  \le 4/  b^{i-1} $.

For instance, the probability of a number $n < 10^9$ being
9-palindromic in base 10 is $0.00009$, while the probability of it
being 9-palindromic in any base is below $0.004$.

\smallskip

From the point of view of probability, the situation is then --- in most
cases --- quite clear:
for $k \le 2 $ every number is $k$-palindromic, while
for $k \ge 4$ almost every number is  not.

\smallskip

The critical case is  $k:=3$.
Consider the following table:

\medskip

\hrulefill

\medskip

\begin{tabular}{ccc}
$\Phi_3(10^2 + 100) - \Phi_3(100) $  \  & = &  \ 61 \\[2mm]
$\Phi_3(10^3+ 100) - \Phi_3(10^3) $  \  & = & \ 70  \\[2mm]
$\Phi_3(10^4+ 100) - \Phi_3(10^4) $  \    & = & \ 83 \\[2mm]
$\Phi_3(10^5+ 100) - \Phi_3(10^5) $  \   & = & \ 86 \\[2mm]
$\Phi_3(10^6+ 100) - \Phi_3(10^6) $  \   & = & \ 89 \\[2mm]
$\Phi_3(10^7+ 100) - \Phi_3(10^7) $   \   &  = & \ 94
\end{tabular}

\medskip

\hrulefill

\medskip

This suggests that almost every number is $3$-palindromic, but { \em not } every
sufficiently large number.
To tackle this problem
it might be worth considering the following reformulation.
We recall that
$\{ \xi \} := \xi - [\xi] \in [0, 1)$ denotes
the fractional part
of a real number $ \xi \in \R$.

\begin{undef}{Lemma 2}
Let $n, b \in \N$ such that
the $b$-expansion of $n$ has length 3.
Then $n$ is 3-palindromic in base $b$ if and only if
$$
\{\,  (n+1) \, \frac{b }{ b^2 + 1} \, \} < \frac{b}{b^2 + 1}.
$$

\end{undef}

\vspace{-4mm}

\begin{proof}{Proof.--}

First note that
the hypothesis
that the $b$-expansion of $n$ has length 3
is equivalent to the fact that
$b^2 +1\le n \le b^3-1$.
Now $n$ is $3$-palindromic in base $b$ if and only if there exists
$0 <  e < b $ and  $0 \le f < b$ such that
$$
n = e \, ( b^2+1) + f \, b .
$$
Solving the associated Diophantine
linear equation $n = x \, (b^2+1) + y \, b$ with respect to
 $x,y$
we see that the above representation is equivalent to the
existence of $ \ell \in \Z$ satisfying
$$
0 < n - \ell \,  b < b  \quad \quad , \quad \quad
0 \le
\ell \, (b^2+1) - n \, b < b.
$$
The second pair of inequalities is equivalent to
$ n \, b / (b^2 + 1) \le \ell < (n+1) \, b / (b^2 + 1)$, and
so it
implies that
$\{ \, (n+1) \, b / (b^2 + 1) \, \} < b / (b^2 + 1)$.
Then this condition is necessary for $n$ to be
3-palindromic in base $b$.

\smallskip

Let's check that it is also sufficient:
the integer $\ell:= [(n+1) \, b / (b^2+1)]$ satisfies the second
pair of inequalities.
Then it
only remains to prove
that it also satisfies the first pair, which is
equivalent to $\ell < n/b < \ell+1$.
This follows from the inequalities
$$
(n+1) \, \frac{b}{b^2+1} < \frac{n}{b} < n \,\frac{ b}{b^2+1} + 1,
$$
which are in term a consequence of the hypothesis
$b^2 +1\le n \le b^3-1$.

\end{proof}

\bigskip

Now we  begin to look at palindromicity as an intrinsic property
--- not attached to any particular base ---
nothing stop us from considering the fact that a given number
can be palindromic in several different basis.
For instance
$$
3074 = (44244)_5 = (22122)_6.
$$

Common sense dictates that multiple palindromicity should be a much more rare
occurrence than simple one, which is also rare as we have already
shown.
In fact it even seems unclear whether
there are numbers which are
$k$-palindromic in as many basis as desired.
We formalize this: let
$$ \mu_k(n) := \# \{ \ b  \, ; \, n
\mbox{ is } k \mbox{-palindromic in base } b \ \} .
$$

So in first instance, we propose the problem of determining whether  $\mu_k$ is
unbounded or not.
Again the cases $k=1,2$ are easy. In the first case
$n=(n)_m$ for any base $m > n$, and so $\mu_1(n) = \infty$ for
every
$n$.
In the second case, set $n:= 2^{2 \,u+1} $ for some $u \in \N$.
Then $n = 2^v \, (2^w-1) + 2^v = (2^v , 2^v)_{2^w-1} $
for $v < w$ such that  $v+w = 2 \, u+1 $.
Then $\mu_2(n) \ge u$.

\smallskip

The following solves the case $k=3$:

\begin{undef}{Theorem 3}
There exists an infinite  sequence $n_1 <n_2 < n_3 < \cdots$ such that
$$
\mu_3(n_j) \ge \frac{1}{7} \, \log (n_j+1).
$$
\end{undef}

\vspace{-10mm}

\begin{proof}{Proof.--}

Take $N \gg 0 $ and assume that
$\mu_3(n) < (1/7) \, \log (N+1)$ for all $n \le N$, so that
\begin{equation} \label{N/7}
\sum_{b} \Phi_3(N, b)
= \sum_{n=1}^N \mu_3(n)
\, <
\frac{1}{7}\, N \, \log (N+1).
\end{equation}
We will see in a minute that this is contradictory:

Set  $ L:= [(N+1)^\frac{1}{3}]$ and $M:= [(N-1)^\frac{1}{2} ]$.
For $ L\le b \le M$
we let $\zeta (b) \in \N$ be the largest integer such that
$$
\zeta(b)  \, (b^2+1) + (b-1) \, b \le N.
$$
Then $1 \le \zeta (b) \le b-1$, and also
every 3-palindromic number $ n:= e  \, b^2 + f \, b + e$  with
$e \le \zeta(b)$ is less or equal than $N$.
Hence
$ \Phi_3(N, b)
 \ge
\zeta(b) \, b$ which implies that
$$
\sum_b \Phi_3(N, b)
\, \ge  \,
\sum_{b=L}^M \Phi_3(N, b)
\, \ge \,
\sum_{b=L}^M  \zeta(b) \, b
\, \ge \,
\sum_{b=L}^M  b\, (  \frac{N}{b^2+1} -2),
$$
as $\zeta(b) +2 \ge N / (b^2+1) $.
We have that
$$
\sum_{b=L}^M  b\, (  \frac{N}{b^2+1}
-2)
\ge N \, \int_{L}^{M+1} \frac{t}{t^2+1} \, dt - 2 \, M^2 \ge
\frac{N}{2} \,
\left(
\log ((M+1)^2 +1) - \log (L^2+1)
\right)
- 2\, N .
$$

We have that $(M+1)^2 +1 \ge  N$ and $L^2 +1 \le 2 \,
N^{2/3}$ and thus we conclude $
\sum_p \Phi_3(N, p) > \frac{N}{6} \, \log N - 2 \, N - \log 2$,
which contradicts Inequality \ref{N/7} for  $N$ large enough.

\smallskip

It follows that for each (sufficiently large) $N \in \N$
there exists $n \le
N$ such that
$$
\mu_3(n)
\, \ge  \,
\frac{1}{7} \, \log (N+1)
\, \ge \,
\frac{1}{7}\, \log (n+1).
$$
The fact that
$\mu_3(n) <\infty$ implies that the set of such $n$'s is infinite.

\end{proof}

Here is some sample data for the cases $k:=4, 5$:
$$
\mu_4(624) = \mu_4(910) = 2
\quad \quad , \quad \quad
\mu_4(19040) = 3
\quad \quad , \quad \quad
\mu_5(2293)= 2.
$$

For $k, \ell , N\in \N$, we
let $\Phi_{k, \ell} (N)$ be the number of $n \le N$ which are
$k$-palindromes
in $\ell$ different basis, that is
$$
\Phi_{k, \ell} (N) := \# \{
\ n \le N \, ; \,  \mu_k (n) \ge \ell\ \ \}.
$$
In particular $\Phi_{k, 1} = \Phi_k$.
The following table gives some more informative data:

\medskip

\begin{center}

\begin{tabular}{|c|c|c|c|}
\hline
&&&\\[-3.7mm]
\hspace*{10mm} $k$ \hspace*{10mm} &
\hspace*{10mm} $\ell$  \hspace*{10mm} &
\hspace*{10mm}
$N$
\hspace*{10mm} &
\hspace*{15mm} $ \Phi_{k, \ell} (N)  \hspace*{15mm}$
\\[0.3mm]
\hline \hline
&&& \\[-3.7mm]
4 & 2 & $10^4 $& 13 \\[0.3mm]
\hline
&&& \\[-3.7mm]
4 & 3 & $10^5 $ & 2 \\[0.3mm]
\hline
&&& \\[-3.7mm]
4 & 4&  $10^5 $ & 0 \\[0.3mm]
\hline
&&& \\[-3.7mm]
5 & 2 & $ 10^4 $ & 10
\\[0.3mm]
\hline
&&& \\[-3.7mm]
5 &  3 & $10^5 $ & 0 \\[0.3mm]
\hline
&&& \\[-3.7mm]
6 & 2 & $ 10^5
$ & 0 \\[0.3mm]
\hline
\end{tabular}

\end{center}

\bigskip

This suggest that
for $k\ge 4$ and  $k + \ell \ge 8$
there are no $k$-palindromic numbers with multiplicity
$\ell$ at all.

\medskip

Finally we can also consider
$$
\mu_{\ge k}(n) := \# \{ \ b  \, ; \, n
\mbox{ is } j \mbox{-palindromic in base } b \mbox{ for some } j \ge k \ \}
= \sum_{j \ge k} \mu_j(n),
$$
that is the number of different basis in which $n$ is
a palindrome of length {\em at least  } $k$.
It is easy to see that this function is unbounded:
we have that
\vspace{-3mm}
$$
n_L:= 2^{2^L}-1 =
( \overbrace{2^{2^\ell}-1, \cdots , 2^{2^\ell}-1 }^{2^{L-\ell}})_{2^{2^\ell}}
$$
and so $n_L$ is $2^{L-\ell}$-palindromic in base $2^{2^\ell} $ for $\ell=0,
\dots, L$.
Hence
$\mu_{\ge k} (n_L) \ge L - \log_2  k $.

\bigskip

A further problem is to determine the {\em density } of
$k$-palindromic numbers in $\ell$ different basis.
From this point of view, Theorem 1 is an important
advance towards the solution of the cases $k \ge 4 $ and $\ell=1$.

The cases when  $\ell \ge 2$ seem to be   much more elusive, but
also interesting.
A solution of them would allow you, for instance, to know how lucky you are
when the number of the taxi-cab you are riding is the
$$
19040 = (8888)_{13} = (5995)_{15} = ( 2 , 14,  14,  2)_{19} \, .
$$

\typeout{Referencias}

\bigskip

\noindent {\sc Antonio J. Di Scala: }
Facultad de Matem{\'a}tica, Astronom{\'\i}a y F{\'\i}sica
(Fa.M.A.F.), Universidad Nacional de C{\'o}rdoba, Ciudad
Universitaria, 5000 C{\'o}rdoba, Argentina \\[1.5mm]
{ \tt E-mail: discala@mate.uncor.edu }

\vspace{4
mm}

\noindent {\sc Mart{\'\i}n Sombra: }
Universit{\'e} de Paris 7, UFR de Math{\'e}matiques,
{\'E}quipe de G{\'e}om{\'e}trie et Dinamique, 2 place Jussieu, 75251 Paris Cedex 05,
France; and 
Departamento de Matem{\'a}tica,
Universidad Nacional de La Plata,
Calle 50 y 115,
1900 La Plata, Argentina. \\[1.5mm]
{\tt E-mail: sombra@jussieu.math.fr}

\end{document}